\documentclass[a4paper]{article}

\usepackage[latin1]{inputenc}
\usepackage[francais]{babel}
\usepackage{delarray,verbatim,enumerate}
\usepackage{amsmath,amsthm,amstext,amsbsy,amssymb,amsfonts,amscd,bm}

\usepackage{ifpdf}
\ifpdf
	\usepackage{aeguill}
\else
	\usepackage[T1]{fontenc}
\fi

\newtheorem{Th}{Th\'eor\`eme}[]

\newtheorem{Cor}[Th]{Corollaire}

\newtheorem{Def} [Th]{D\'efinition}

\newtheorem{ThD}[Th]{Th\'eor\`eme \& D\'efinition}

\def\PreuveTh{\smallskip\noindent {\it Preuve du Théorème.~}}
\def\Remarque{\smallskip\noindent {\it Remarque.~}}

\font\teneufm=eufm10
\font\seveneufm=eufm7
\font\fiveeufm=eufm5
\newfam\gothfam
\textfont\gothfam=\teneufm \scriptfont\gothfam=\seveneufm
\scriptscriptfont\gothfam=\fiveeufm
\def\goth{\fam\gothfam}

\def\QQ{\mathbb Q}	
			\def\ZZ{\mathbb Z}	
\def\Zl{\mathbb Z_\ell}		\def\Fl{\mathbb F_\ell}
\def\F2{\mathbb F_2}		\def\Z2{\mathbb{Z}_2}	
\def\Z2{\mathbb{Z}_2}

  		\def\C{\mathcal  C}		
 	  	\def\Cl{\mathcal  C \ell}
		
\def\G{\mathcal  G}		\def\F{\mathcal  F}	\def\T{\mathcal T}

		\def\p{{\goth p}}		\def\l{{\goth l}}

\def\rg{\operatorname{rg}}	
	
\def\Gal{\operatorname{Gal}}

\newcommand{\cdast}{\mbox{\large $\circledast$}}
\def\oast{\operatorname{\circledast}}
\def\Oast{\operatornamewithlimits{\cdast}}
\newcommand{\cdplus}{\mbox{\large $\oplus$}}
\def\Oplus{\operatornamewithlimits{\cdplus}}

\begin{document}

\title{\LARGE\bf Note sur les corps 2-rationnels}

\author{ Jean-François {\sc Jaulent} }
\date{}
\maketitle
\bigskip

{\small
\noindent{\bf Résumé.} Nous déterminons le groupe de Galois de la pro-2-extension 2-ramifiée maximale d'un corps de nombres 2-rationnel.
}

\

{\small
\noindent{\bf Abstract.} We compute the Galois group of the maximal 2-ramified pro-2-extension of a 2-rational number field.
}
\bigskip\bigskip

\bigskip
%%%%%%%%%%%%%%%%%%%%%%%%%%%%%%%%%%%%%%%%%%%%%%%%%%%%%%%%

\noindent{\large \bf Introduction}

%%%%%%%%%%%%%%%%%%%%%%%%%%%%%%%%%%%%%%%%%%%%%%%%%%%%%%%%
\medskip

Les notions de corps $\ell$-rationnel ou $\ell$-régulier (pour un nombre premier $\ell$)  introduites idépendamment par A. Movahhedi et T. Nguyen Quang Do dans \cite{MN} d'une part, par G. Gras et l'auteur dans \cite{GJ} d'autre part, se trouvent coïncider en présence des racines $\ell$-ièmes de l'unité\footnote{En fait dès que le corps $F$ considéré contient le sous-corps réel $\QQ[\zeta_\ell+\bar\zeta_\ell]$ du corps $\QQ[\zeta_\ell]$.} donc, tout spécialement, pour $\ell=2$.\smallskip

\begin{itemize}

\item La $\ell$-régularité se définit naturellement en termes de $K$-théorie et exprime simplement la trivialité du $\ell$-noyau régulier du corps $F$ considéré ({\em i.e.} du noyau dans la $\ell$-partie du groupe universel $K_2(F)$ des symboles de Hilbert attachés aux places non complexes qui ne divisent pas $\ell$).\smallskip

\item La $\ell$-rationalité s'exprime en termes galoisiens et traduit la pro-$\ell$-liberté du groupe de Galois $\G_F=\Gal(M_F/F)$ attaché à la pro-$\ell$-extension (galoisienne) $\ell$-ramifiée $\infty$-décomposée maximale $M_F$ de $F$ ({\em i.e.} au compositum des $\ell$-extensions de $F$ qui sont non ramifiées aux places finies\footnote{En accord avec les conventions de la Théorie $\ell$-adique du Corps de Classes ({\em cf.} \cite{G2,Ja}), nous ne parlons jamais de ramification à l'infini, mais de {\em complexification} des places réelles.} étrangères à $\ell$ et complètement décomposée aux places à l'infini).
\end{itemize}\smallskip

Plus précisément, d'après \cite{JN} (Th.1.2), la $\ell$-rationalité d'un corps de nombres $K$ s'exprime comme suit :

\setcounter{Th}{-1}
\begin{ThD} Les conditions suivantes sont équivalentes :
\begin{itemize}
\item[(i)] Le groupe de Galois $\mathcal G_K = \Gal (M_K/K)$ de la pro-$\ell$-extension $\ell$-ramifiée $\infty$-décomposée maximale $M_K$ de $K$ est un pro-$\ell$-groupe libre (sur $d_K=1 + c_K$ générateurs, si $c_K$ désigne le nombre de places complexes de $K$).

\item[(ii)]  Le groupe de Galois $\G^{ab}_K = \Gal (M_K^{ab}/K)$ de la sous-extension  abélienne maximale $M_K^{ab}/K$ de $M_K/K$  est un $\Zl$-module libre de dimension $1+c_K$.

\item[(iii)]  Le corps $K$ v\'erifie la conjecture de Leopoldt (pour le premier $\ell$) et le sous-module de torsion $\mathcal T_K$ du groupe $\G^{ab}_K = \Gal (M_K^{ab}/K)$ est trivial.

\item[(iv)]  Le groupe $V_K = \{ x \in K^{\times} |   \ x \in K^{\times \ell}_\l \  \forall  \l \mid \ell \ \& \  v_{\goth p} (x) \equiv 0  \mod \ell  \  \ \forall {\p} \nmid \ell\infty \}$\linebreak des éléments $\ell$-hyperprimaires du corps $K$ se réduit à $K^{\times \ell}$, et l'on a l'identité entre les $\ell$-rangs des $\ell$-groupes de racines de l'unité : 

\centerline {$\rg_\ell\, \bm\mu_K = \sum_{{\l} | \ell} \rg_\ell\, \bm\mu_{K_\l}.$}

\item[(v)]  Le corps $K$ v\'erifie l'une des deux conditions suivantes :
\begin{itemize}
\item[(a)]  ou bien $K$ contient une racine primitive $\ell$-ième de l'unité $\zeta$, auquel cas $K$ possède une unique place $\l$ au dessus de $\ell$, et le $\ell$-groupe des $\ell$-classes d'idéaux au sens restreint $\Cl'_K$ est trivial ;

\item[(b)]  ou bien $K$ ne contient pas $\zeta$, auquel cas les places de $K$ au-dessus de $\ell$ ne se décomposent pas complètement dans l'extension cyclotomique $K [\zeta ] /K$ et la $\omega$-composante du $\ell$-groupe des $\ell$-classes
d'idéaux au sens restreint $\Cl'_{K [\zeta ]}$ du corps $K [\zeta]$ est triviale, si $\omega$ désigne le caractère cyclotomique ({\em i.e.} le caractère de l'action sur $\bm\mu_K$ de $\Gal (K [\zeta]/K)$).
\end{itemize}
\end{itemize}
Lorsque ces conditions sont réunies, le corps $K$ est dit $\ell$-rationnel.
\end{ThD}

\Remarque Pour $\ell=2$, il résulte clairement de la condition {\em (v,a)} ci-dessus qu'un corps 2-rationnel ne peut contenir qu'{\em une seule} place au-dessus de 2.\medskip

De fait les prémices de la notion de $\ell$-régularité remontent aux travaux de G. Gras notamment à sa note sur le $K_2$ des corps de nombres \cite{G1}, tandis que la notion de $\ell$-rationalité apparait (sous une forme cachée) dans les travaux de H. Miki \cite{Mi}, à l'occasion de l'étude d'une condition suffisante de la conjecture de Leopoldt, ainsi que dans ceux de K. Wingberg \cite{W1,W2}, où est étudiée la même condition.\smallskip

Les articles \cite{GJ,MN} cités plus haut en caractérisent complètement la propagation (dans une $\ell$-extension) en termes de primitivité de la ramification, ce que les approches antérieures ne donnaient pas. Une synthèse de leurs résultats est présentée dans \cite{JN} et un exposé systématique en est donné dans le livre de G. Gras \cite{G2}. \smallskip

Diverses généralisations de ces notions ont été étudiées par O. Sauzet et l'auteur ({\em cf}. \cite{JS1,JS2}), notamment dans le cas $\ell=2$ qui se révèle, comme à l'ordinaire, le plus compliqué ; elles donnent naissance, en particulier, à la notion de corps 2-birationnel. \par
Tout récemment, J. Jossey \cite{Jo} a jugé bon d'introduire une notion de corps $\ell$-rationnel qui diffère de celles déjà utilisée (en fait pour $\ell=2$, dès lors que le corps de nombres considéré posssède des plongements réels), ce qui parait doublement malheureux, car s'écartant de l'usage établi et ne s'appliquant pas au corps {\em des} rationnels $\QQ$.\par
Pour ces raisons, afin de prévenir toute confusion, il nous semble plus judicieux de parler dans son contexte de corps {\em 2-surrationnel}. Précisons ce point :

\begin{Def}
Soient $K$ un corps de nombres ayant $r_K$ places réelles et $c_K$ places complexes, $\ell$ un nombre premier, $M'_K$ la pro-$\ell$-extension  maximale de $K$ qui est 2-ramifiée ({\em i.e.} non ramifiée aux places (finies) étrangères à $\ell$) et $M_K$ la sous-extension maximale de $M'_K$ qui est complètement décomposée aux places à l'infini. Nous disons que le corps $K$ est :
\begin{itemize}
\item[(i)] $\ell$-surrationnel, lorsque le groupe $\G'_K=\Gal(M'_K/K)$ est pro-$\ell$-libre ;
\item[(ii)] $\ell$-rationnel, lorsque son quotient $\G_K=\Gal(M_K/K)$ est pro-$\ell$-libre.
\end{itemize}
\end{Def}

On voit alors que les deux notions coïncident, sauf précisément dans le cas où $\ell$ vaut 2 et où le corps $K$ possède une ou plusieurs places réelles se complexifiant dans $M'_K$ : ainsi tout  corps $\ell$-surrationnel est-il évidemment $\ell$-rationnel, mais la réciproque n'est pas vraie en présence de complexification.
\bigskip

Le but de cette note est de déterminer la structure du groupe $\G'_K$ dans le cas exceptionnel où le corps considéré $K$ est 2-rationnel mais non 2-surrationnel.

%%%%%%%%%%%%%%%%%%%%%%%%%%%%%%%%%%%%%%%%%%%%%%%%%%%%%%%%%%%
\bigskip\medskip

\noindent{\large \bf Théorème principal : description du groupe $\G'_K$}

\bigskip
%%%%%%%%%%%%%%%%%%%%%%%%%%%%%%%%%%%%%%%%%%%%%%%%%%%%%%%%%%%

Notre résultat est extrèmement simple, qui s'énonce comme suit :

\begin{Th}
Soit $K$ un corps de nombres 2-rationnel possédant $r_K \ge 1$ places réelles et $c_K$ places complexes. Le groupe de Galois $\G'_K=\Gal(M'_K/K)$ de la pro-2-extension 2-ramifiée maximale $M'_K$ de $K$ est alors le pro-2-produit libre
$$
\G'_K \; \simeq \; \Z2^{\oast (1+c_K)} \;\Oast\; C_2^{\oast (r_K)}
$$
de $1+c_K$ exemplaires du groupe procyclique $\Z2$ et de $r_K$ exemplaires du groupe cyclique $C_2\simeq\ZZ/2\ZZ$.
\end{Th}

\begin{Cor}
Les corps de nombres 2-rationnels qui sont 2-surrationnels sont ceux totalement imaginaires.
\end{Cor}

\PreuveTh
Si le corps 2-rationnel $K$ considéré n'admet pas de plongement réel, {\em i.e.} dans le cas $r_K=0$, il est alors 2-surrationnel ; et le groupe $\G'_K=\G_K$ est alors pro-2-libre sur $d_K=c_K+1$ générateurs, comme annoncé.
\smallskip

Sinon, {\em i.e.} dans le cas $r_K>0$, introduisons l'extension quadratique $L=K[i]$ engendrée par les racines 4\iemes {} de
l'unité. Elle est évidemment 2-ramifiée sur $K$ donc, en vertu des théorèmes de propagation de \cite{GJ,MN} ({\em cf. e.g.} \cite{JN}, Th. 3.5), 2-rationnelle ; et finalement 2-surrationnelle, puisque totalement imaginaire. En d'autres termes, le groupe de Galois $\G_L=\G'_L$ de la pro-2-extension 2-ramifiée maximale $M_L$ de $L$ est pro-2-libre :
$$
\G_L \; \simeq  \; \Z2^{\oast d_L}, \quad {\rm avec } \quad d_L = c_L +1 = r_K + 2c_K +1.
$$

Cela étant, comme l'extension quadratique $L/K$ est 2-ramifiée, $M_L$ est aussi la plus grande pro-2-extension $M'_K$ de $K$ qui est 2-ramifiée et le groupe de Galois $\G'_K$ est ainsi {\em potentiellement libre}, puisqu'il contient un sous-groupe ouvert qui est pro-2-libre, à savoir $\G_L$, lequel est d'indice 2 dans $\G'_K$.\par
Plus précisément, les résultats de structure de W. Herfort et P. Zalesskii ({\em cf. } \cite{HZ}, Th. 0.2) assurent qu'il existe alors une famille finie $(\F_i)_{i=0,\cdots ,k}$ de pro-2-groupes libres sur respectivement $d_0,\cdots, d_k$ générateurs, où $k$ dénombre les classes de conjugaisons de sous-groupes d'ordre 2 de $\G'_K$, tels qu'on ait :
$$
\G'_K \; \simeq \; \F_0 \Oast \left( \Oast_{i=1}^k (C_2 \times \F_i) \right) .
$$
En particulier, l'abélianisé $\G^{'ab}_K$ de $\G'_K$ admet alors la décomposition directe :\smallskip

\centerline{$\G^{'ab}_K \; \simeq \; \Z2^{d_0} \oplus \left( \Oplus_{i=1}^k (C_2 \oplus \Z2^{d_i}) \right) \; \simeq \; \Z2^{\sum d_i} \oplus C_2^k$.}\smallskip

Et, puisque le corps 2-rationnel $K$ vérifie la conjecture de Leopoldt, il vient :\medskip

\centerline{$\sum_{i=0}^k d_i = d_K = c_K+1$ ; ainsi que l'isomorphisme : $\T'_K \simeq \C_2^k$.}\medskip

\noindent où $\T '_K$ désigne le sous-groupe de $\Z2$-torsion du groupe $\G^{'ab}_K$. \smallskip

Maintenant, il résulte clairement de la pro-2-décomposition de $\G'_K$ que les nombres minimaux de générateurs $d(\G'_K)$ et de relations $r(\G'_K)$ qui le définissent comme pro-2-groupe sont respectivement :\smallskip

\centerline{$d(\G'_K) = k + \sum_{i=0}^k d_i = k+c_K  +1\;$ et $\;r(\G'_K) = \sum_{i=1}^k (1+d_i) = d(\G'_K) - d_0$.}\pagebreak

\noindent Or les formules de \v Safarevi\v c ({\em cf.} \cite{Sa}, \cite{G2}, § III.4, ou encore \cite{NW}, Th. 8.7.3) :
$$
\aligned
d ( \G'_K) &= \dim_{\Fl} (H^1 (\G'_K, \Fl))\\
&= c_K + 1 + \dim_{\Fl} V_K/K^{\times \ell} + ( \sum\limits_{{\l}|\ell} \rg_\ell\,\bm\mu_{K_\l} - \rg_\ell\,\bm\mu_K)~;~ {\rm et}\\
r(\G'_K) &= \dim_{\Fl} (H^2(\G'_K, \Fl)) = \dim_{\Fl}V_K/K^{\times \ell} + ( \sum_{\l | \ell} \rg_\ell\,\bm\mu_{K_\l} - \rg_\ell\,\bm\mu_K).
\endaligned
$$ 
donnent ici, d'après la condition $(v,a)$, pour $\ell=2$ et $K$ supposé 2-rationnel :
$$
d(\G'_K) = \dim_{\mathbb F_2}{}^2\G_K^{'ab} = r_K + c_K + 1 \qquad {\rm et} \qquad r(\G'_K) = \dim_{\mathbb F_2}{}^2\T'_K =  r_K\; ;
$$
d'où finalement : $k=r_K$ et $d_0=1+c_K$ ; ce qui conduit au résultat attendu.

%%%%%%%%%%%%%%%%%%%%%%%%%%%%%%%%%%%%%%%%
%REFERENCES
%%%%%%%%%%%%%%%%%%%%%%%%%%%%%%%%%%%%%%%%

\def\refname{\small{\sc  Références}}

\bigskip\noindent
{\small
\begin{tabular}{l}
{Jean-Fran\c cois {\sc Jaulent}}\\

Universit{\'e} de Bordeaux\\
Institut de Math{\'e}matiques de Bordeaux \\
351, cours de la lib{\'e}ration\\
F-33405 {\sc Talence} Cedex\\
courriel : Jean-Francois.Jaulent@math.u-bordeaux1.fr 
\end{tabular}
}

 \end{document}